\documentclass[final,twoside,11pt]{entics}
\usepackage{enticsmacro}
\usepackage{graphicx}
\usepackage[all]{xy}
\usepackage{hyperref}
\usepackage{dsfont}
\usepackage{tikz}
\usepackage{float}
\usepackage{bm}
\usepackage{mathrsfs}
\usepackage{amsmath}
\newcommand{\ua}{\mathord{\uparrow}}
\newcommand{\da}{\mathord{\downarrow}}
\newcommand{\rom}[1]{\rm{\uppercase\expandafter{\romannumeral #1}}}

\def\conf{ISDT 2022} 	%%Fill in the acronym for your conference (with year)
\volume{2}			%Fill in the ENTICS volume number here
\def\volu{2}			% and here
\def\papno{5}			%Fill in your paper number here
\def\mydoi{{\normalshape  \href{https://doi.org/10.46298/entics.10349}{10.46298/entics.10349}}}
\def\copynames{M. Jin, Q. Li} %% Fill in the first initial and last name of the authors
\def\runtitle{On $k$-ranks of topological spaces}

\def\lastname{Jin and Li}
\begin{document}
\begin{frontmatter}
  \title{On $k$-ranks of Topological Spaces}
  \author{Mengjie Jin\thanksref{mengjiejinjin@163.com}}	
   \author{Qingguo Li\thanksref{ALL}\thanksref{liqingguoli@aliyun.com}}	
   \address{School of Mathematics\\ Hunan University\\			
    Changsha, Hunan, 410082, China}  							
  \thanks[ALL]{This work is supported by the National Natural Science Foundation of China (No.12231007)}   
  \thanks[mengjiejinjin@163.com]{Email: \href{mailto:mengjiejinjin@163.com} {\texttt{\normalshape
        mengjiejinjin@163.com}}}
  \thanks[liqingguoli@aliyun.com]{Corresponding author, Email:  \href{mailto:liqingguoli@aliyun.com} {\texttt{\normalshape
        liqingguoli@aliyun.com}}}
\begin{abstract}
  In this paper, the concepts of $K$-subset systems and $k$-well-filtered spaces are introduced, which provide another uniform approach to $d$-spaces, $s$-well-filtered spaces (i.e., $\mathcal{U}_{S}$-admissibility) and well-filtered spaces. We prove that the $k$-well-filtered reflection of any $T_{0}$ space exists. Meanwhile, we propose the definition of $k$-rank, which is an ordinal that measures how many steps from  a $T_{0}$ space to a $k$-well-filtered space. Moreover, we derive that for any ordinal $\alpha$, there exists a $T_{0}$ space whose $k$-rank equals to $\alpha$. One immediate corollary is that for any ordinal $\alpha$, there exists a $T_{0}$ space whose $d$-rank (respectively, $wf$-rank) equals to $\alpha$.
\end{abstract}
\begin{keyword}
  $k$-well-filtered space, $k$-Rudin set, $k$-well-filtered reflection, $k$-rank
\end{keyword}
\end{frontmatter}
\section{Introduction}\label{intro}
In non-Hausdorff topological spaces and domain theory, $d$-spaces and well-filtered spaces are two important classes of spaces. Let $\mathbf{Top}_{0}$ be the category of all $T_{0}$ spaces, $\mathbf{Top}_{d}$ the category of all $d$-spaces and $\mathbf{Top}_{w}$ the category of all well-filtered spaces. It is well-known that $\mathbf{Top}_{d}$ and $\mathbf{Top}_{w}$ are reflective in $\mathbf{Top}_{0}$, respectively. Different ways for constructing $d$-completions and well-filtered reflections of $T_{0}$ spaces were found in \cite{Ershov 17,Liu2020,Shen19,Zhongxi}. In \cite{Ershov 17}, Ershov introduced one way to get $d$-completions of $T_{0}$ spaces using the equivalent classes of directed subsets, he called it $d$-rank which is an ordinal that measures how many steps from  a $T_{0}$ space to a $d$-space. Inspired by his method, in \cite{Liu2020}, Liu, Li and Wu proposed one way to get well-filtered reflections of $T_{0}$ spaces using the equivalent classes of Rudin subsets, they called it $wf$-rank, which is an ordinal that measures how far a $T_{0}$ space is from being a well-filtered space.

In \cite{Xiaoquan21}, based on irreducible subset systems, Xu provided a uniform approach to $d$-spaces, sober spaces and well-filtered spaces, and developed a general framework for dealing with all these spaces. In this paper, we will provide another uniform approach to $d$-spaces and well-filtered spaces and develop a general framework for dealing with all these spaces. Similar to the concept of irreducible subset systems in~\cite{Xiaoquan21}, we propose the concepts of $K$-subset systems and $k$-well-filtered spaces. For a $K$-subset system $Q_{k}:\mathbf{Top}_{0}\longrightarrow \mathbf{Set}$ and a $T_{0}$ space $X$, $X$ is called $k$-well-filtered if for any open set $U$ and a filtered family $\mathcal{K}\subseteq Q_{k}(X)$, $\bigcap\mathcal{K}\subseteq U$ implies $K\subseteq U$ for some $K\in \mathcal{K}$. The category of all $k$-well-filtered spaces with continuous mappings is denoted by $\mathbf{Top}_{k}$. It is not difficult to verify that $d$-spaces and well-filtered spaces are two special kinds of $k$-well-filtered spaces. Moreover, we find that $s$-well-filtered spaces (i.e., $\mathcal{U}_{S}$-admissibility in \cite{Heckmann92}) is also a kind of $k$-well-filtered spaces which is different from $d$-spaces and well-filtered spaces. Just like directed subsets and Rudin subsets, the concept of $k$-Rudin sets will be introduced. Moreover, for a $K$-subset system $Q_{k}:\mathbf{Top}_{0}\longrightarrow \mathbf{Set}$, we use the equivalent classes of $k$-Rudin sets to construct the $k$-well-filtered reflections of $T_{0}$ spaces. Meanwhile, we introduce the concept of $k$-rank, which is an ordinal that measures how far a $T_{0}$ space can become a $k$-well-filtered space. For a $T_{0}$ space $X$, we get that there exists an ordinal $\alpha$ such that the $k$-rank of $X$ is equal to $\alpha$.

In \cite{Ershov 17} and \cite{Liu21}, for any ordinal $\alpha$, there exists a $T_{0}$ space whose $d$-rank (respectively, $wf$-rank) equals to $\alpha$. Consider a $T_{0}$ space whose $k$-rank equals to $\alpha$ may be more complex, because we know little about $Q_{k}(X)$. We have to find suitable conditions to characterize a class of $T_{0}$ spaces whose $k$-rank equals to $\alpha$. It turns out that finding these $T_{0}$ spaces is the hard part of our task, but how to prove the results is relatively simple.

Finally, we obtain that for any ordinal $\alpha$, there exists a $T_{0}$ space whose $k$-rank equals to $\alpha$.

\section{Preliminaries}
First, we briefly recall some standard definitions and notations to be used in this paper, for further details see \cite{Davey02}, \cite{gierz03}, \cite{nht2} and \cite{Heckmann13}.

Let $P$ be a poset and $A\subseteq P$. We denote $\ua A=\{x\in P \mid x\geq a \mbox{ for some } a\in A\}$ and $\da A=\{x\in P \mid x\leq a \mbox{ for some } a\in A\}$. For every $a\in P$, we denote $\ua\{a\}=\ua a=\{x\in P \mid x\geq a\}$ and $\da\{a\}=\da a=\{x\in P \mid x\leq a\}$. $A$ is called an \emph{upper set} (resp., a \emph{lower set}) if $A=\ua A$ (resp., $A=\da A$). $A$ is called \emph{directed} provided that it is nonempty and every finite subset of $A$ has an upper bound in $A$. The set of all directed sets of $P$ is denoted by $\mathcal{D}(P)$. Moreover, the set of all nonempty finite sets in $P$ is denoted by $P^{<\omega}$.

A poset $P$ is called a \emph{dcpo} if every directed subset $D$ in $P$ has a supremum. A subset $U$ of $P$ is called \emph{Scott open} if (1) $U=\ua U$ and (2) for any directed subset $D$ for which $\vee D$ exists, $\vee D\in U$ implies $D\cap U\neq \emptyset$. All Scott open subsets of $P$ form a topology, we call it the \emph{Scott topology} on $P$ and denoted by $\sigma(P)$.

For a $T_{0}$ space $X$, let $\mathcal{O}(X)$ (resp., $\Gamma(X)$) be the set of all open subsets (resp., closed subsets) of $X$. For a subset $A$ of $X$, the closure of $A$ is denoted by $\mathrm{cl}(A)$ or $\overline{A}$. We use $\leq_{X}$ to represent the specialization order of $X$, that is, $x\leq_{X}y$ iff $x\in \overline{\{y\}}$. A subset $B$ of $X$ is called \emph{saturated} if $B$ equals the intersection of all open sets containing it (equivalently, $B$ is an upper set in the specialization order). Let $S(X)=\{\{x\}\mid x\in X\}$, $S_{c}(X)=\{\da x\mid x\in X\}$ and $\mathcal{D}_{c}(X)=\{\overline{D}\mid D\in \mathcal{D}(X)\}$. A $T_{0}$ space $X$ is called a \emph{$d$-space} (i.e., \emph{monotone convergence space}) if $X$ (with the specialization order) is a $dcpo$ and $\mathcal{O}(X)\subseteq \sigma(X)$ (\cite{gierz03}). The category of all $d$-spaces with continuous mappings is denoted by $\mathbf{Top}_{d}$.

For a $T_{0}$ space $X$, let $\mathcal{K}$ be a filtered family under the inclusion order in $Q(X)$, which is denoted by $\mathcal{K}\subseteq_{filt}Q(X)$, i.e., for any $K_{1},K_{2}\in Q(X)$, there exists $K_{3}\in Q(X)$ such that $K_{3}\subseteq K_{1}\cap K_{2}$. $X$ is called \emph{well-filtered} if for any open subset $U$ and any $\mathcal{K}\subseteq_{filt}Q(X)$, $\bigcap \mathcal{K}\subseteq U$ implies $K\subseteq U$ for some $K\in \mathcal{K}$. The category of all well-filtered spaces with continuous mappings is denoted by $\mathbf{Top}_{w}$
(\cite{Xiaoquan20}).

In what follows, $\mathbf{K}$ always refers to a full subcategory $\mathbf{Top}_{0}$ that contains $\mathbf{Sob}$, the full subcategory of sober spaces.  The objects of $\mathbf{K}$ are called \emph{$\mathbf{K}$-spaces}.

\begin{definition}\cite{Xiaoquan20}\label{K-reflction}
Let $X$ be a $T_{0}$ space. A \emph{$\mathbf{K}$-reflection} of $X$ is a pair $\langle \widehat{X},\mu\rangle$ comprising a $\mathbf{K}$-space $\widehat{X}$ and a continuous mapping $\mu$: $X\longrightarrow \widehat{X}$ satisfying that for any continuous mapping $f$: $X\longrightarrow Y$ to a $\mathbf{K}$-space, there exists a unique continuous mapping $f^{\ast}$: $\widehat{X}\longrightarrow Y$ such that $f^{\ast}\circ \mu= f$, that is, the following diagram commutes.
\end{definition}
\centerline{\xymatrix{X \ar[dr]_{f} \ar[r]^{\mu}
                & \widehat{X} \ar[d]^{f^{\ast}}  \\
                &  Y }}

By a standard argument, $\mathbf{K}$-reflections, if they exist, are unique up to homeomorphism. We shall use $X^{k}$ to denote the space of the $\mathbf{K}$-reflection of $X$ if it exists.

For $\mathbf{K}=\mathbf{Top}_{w}$, the $\mathbf{K}$-reflection of X is called the \emph{well-filterification} or \emph{well-filtered reflection} of $X$, we denote it by $\mathrm{H}_{wf}(X)$ if the well-filterification of $X$ exists. For $\mathbf{K}=\mathbf{Top}_{d}$, the $\mathbf{K}$-reflection of X is called the \emph{$d$-completion} of $X$, we denote it by $\mathrm{H}_{d}(X)$ if the $d$-completion of $X$ exists.

\begin{definition}\cite{Xu2016}
Let $X=(X, \tau)$ be a topological space and $A\subseteq X$. $A$ is called strongly compact in $X$ if for each $U\in \tau$ with $A\subseteq U$, there is $F\in X^{<\omega}$ such that $A\subseteq \ua_{\tau}F\subseteq U$.
\end{definition}

\begin{proposition}\cite{Heckmann92}
Every finite set is strongly compact, and every strongly compact set is compact.
\end{proposition}

\begin{proposition}\cite{Heckmann92}
$A$ is strongly compact if and only if $\ua A$ is so.
\end{proposition}

We use $Q_{s}(X)$ to denote the set of all nonempty strongly compact saturated subsets of $X$. $X$ is called \emph{$s$-well-filtered} (i.e., \emph{$\mathcal{U}_{S}$-admissibility} in \cite{Heckmann92}) if it is $T_{0}$, and for any open subset $U$ and $\mathcal{K}\subseteq_{filt} Q_{s}(X)$, $\bigcap\mathcal{K}\subseteq U$ implies $K\subseteq U$ for some $K\in \mathcal{K}$. The category of all $s$-well-filtered spaces with continuous mappings is denoted by $\mathbf{Top}_{s-w}$.

\section{$k$-well-filtered spaces}
In this section, we provide a uniform approach to $d$-spaces and well-filtered spaces and develop a general framework for dealing with all these spaces.

\begin{definition}
$Q_{k}:\mathbf{Top}_{0}\longrightarrow \mathbf{Set}$ is called a \emph{$C$-subset system} if $S^{u}(X)\subseteq Q_{k}(X)\subseteq Q(X)$ for all $X\in ob(\mathbf{Top}_{0})$, where $S^{u}(X)=\{\ua x \mid x\in X\}$.
\end{definition}

\begin{definition}\label{k}
Let $Q_{k}:\mathbf{Top}_{0}\longrightarrow \mathbf{Set}$ be a $C$-subset system and $X$ a $T_{0}$ space. A nonempty subset $A$ is said to have \emph{$k$-Rudin property}, if there exists $\mathcal{K}\subseteq_{filt}Q_{k}(X)$ such that $A$ is a minimal closed set that intersects all members of $\mathcal{K}$. We call such a set \emph{$k$-Rudin} or \emph{$k$-Rudin set}. Let $\mathrm{K}^{R}(X)=\{A\subseteq X \mid A \mbox{ has $k$-Rudin property}\}$ and $\mathrm{K}^{R}_{c}(X)=\mathrm{K}^{R}(X)\cap \Gamma(X)$.
\end{definition}

For $Q_{k}(X)=Q(X)$, a $k$-Rudin set of $X$ is called a \emph{Rudin set} (i.e., KF set) of $X$. The set of all Rudin sets of $X$ is denoted by $\overline{KF}(X)$. $\mathrm{RD}(X)=\overline{KF}(X)\cap \Gamma(X)$.

\begin{proposition}\label{k-Rudin sets}
Let $Q_{k}:\mathbf{Top}_{0}\longrightarrow \mathbf{Set}$ be a $C$-subset system and $X$ a $T_{0}$ space. Then $\mathcal{D}(X)\subseteq \mathrm{K}^{R}(X)\subseteq \overline{KF}(X)$.
\end{proposition}

\begin{proof}
Clearly, $\mathrm{K}^{R}(X)\subseteq \overline{KF}(X)$. Now we prove that every directed subset $D$ of $X$ is a $k$-Rudin set. Let $\mathcal{K}=\{\ua d\mid d\in D\}$. Then $\mathcal{K}\subseteq Q_{k}(X)$ is filtered and $\overline{D}$ interests all members of $\mathcal{K}$. Assume that $A$ is a closed subset in $X$ and interests all members of $\mathcal{K}$. This means that $A\cap \ua d\neq \emptyset$ for all $d\in D$. Since $A$ is closed, it is a lower set, then $d\in A$ for all $d\in D$. Hence, $\overline{D}\subseteq A$. So $\overline{D}$ is a minimal closed set that intersects all members of $\mathcal{K}$.
\end{proof}

\begin{definition}\label{15}
A $C$-subset system $Q_{k}:\mathbf{Top}_{0}\longrightarrow \mathbf{Set}$ is called a \emph{$K$-subset system} provided that for any $T_{0}$ spaces $X$, $Y$ and any continuous mapping $f:X\longrightarrow Y$, $f(A)\in \mathrm{K}^{R}(Y)$ for all $A\in \mathrm{K}^{R}(X)$.
\end{definition}

\begin{definition}\label{k-well-filtered}
Let $Q_{k}:\mathbf{Top}_{0}\longrightarrow \mathbf{Set}$ be a $K$-subset system and $X$ a $T_{0}$ space. $X$ is called \emph{$k$-well-filtered} if for any open set $U$ and $\mathcal{K}\subseteq_{filt}Q_{k}(X)$, $\bigcap\mathcal{K}\subseteq U$ implies $K\subseteq U$ for some $K\in \mathcal{K}$. The category of all $k$-well-filtered spaces with continuous mappings is denoted by $\mathbf{Top}_{k}$.
\end{definition}

In the following, we give some special $k$-well-filtered spaces
and their relations with $\mathbf{Top}_{d}$, $\mathbf{Top}_{w}$ and $\mathbf{Top}_{s-w}$, respectively.

For a $C$-subset system $Q_{k}:\mathbf{Top}_{0}\longrightarrow \mathbf{Set}$ and a $T_{0}$ space $X$, here are some important examples of $Q_{k}(X)$:

\begin{enumerate}
\item [(\texttt{1})] $Q_{k}(X)=S^{u}(X)$ (i.e., $Q_{k}(X)=\{\ua x\mid x\in X\}$).
\item [(\texttt{2})] $Q_{k}(X)=Q_{f}(X)$ (i.e., $Q_{k}(X)=\{\ua F\mid \emptyset\neq F\in X^{<\omega}\}$).
\item [(\texttt{3})] $Q_{k}(X)=Q_{s}(X)$ (i.e., $Q_{k}(X)=\{A\mid A \mbox{ is a nonempty strongly compact saturated subset in } X\}$).
\item [(\texttt{4})] $Q_{k}(X)=Q(X)$ (i.e., $Q_{k}(X)$ is the set of all nonempty compact saturated
subsets in $X$).
\end{enumerate}
Let $Q_{k}:\mathbf{Top}_{0}\longrightarrow \mathbf{Set}$ be a $C$-subset system. For any $T_{0}$ space $X$, if $Q_{k}(X)=S^{u}(X)$ (i.e., $Q_{k}(X)=Q_{f}(X)$), then it follows directly from Definition 1, Example 1(1) and Theorem 1 in \cite{Zhenzhu} that $X$ is $k$-well-filtered iff $X$ is a $d$-space. If $Q_{k}(X)=Q(X)$, it is trivial that $X$ is $k$-well-filtered iff $X$ is well-filtered. In the case $Q_{k}(X)=Q_{s}(X)$, $k$-well-filtered spaces are exactly $s$-well-filtered spaces.

From the above, for a $K$-subset system $Q_{k}:\mathbf{Top}_{0}\longrightarrow \mathbf{Set}$, it is not difficult to see that well-filtered spaces are $k$-well-filtered spaces and $k$-well-filtered spaces are $d$-spaces. That is

$$\mathbf{Top}_{w}\subseteq\mathbf{Top}_{k}\subseteq\mathbf{Top}_{d}.$$

In particular, well-filtered spaces are $s$-well-filtered spaces and $s$-well-filtered spaces are $d$-spaces. In Example \ref{8} and Example \ref{9} below, we will show that the converses are not true, respectively.

\begin{example}\label{8}
Consider set $N$ of natural numbers. Let $X=(N, \tau_{cof})$ be the space $N$ equipped with the co-finite topology (the empty set and the complements of finite subsets of $N$ are open sets). Then
\begin{enumerate}
\item[(\texttt{a})] $\Gamma(X)=\{\emptyset, N\}\cup N^{<\omega}$, $X$ is $T_{1}$, hence it is a $d$-space.
\item[(\texttt{b})] $X$ is $s$-well-filtered since a subset in $T_{1}$ spaces is strongly compact iff it is finite.
\item[(\texttt{c})] $K(X)=2^{N}\backslash \emptyset$.
\item[(\texttt{d})] $\mathrm{RD}(X)=\{N\}\cup \{\{n\}\mid n\in N\}$.
\item[(\texttt{e})] $X$ is not well-filtered.
\end{enumerate}
\end{example}

\begin{example} (Johnstone space)\label{9}
Recall the dcpo constructed by Johnstone in \cite{nht2}, which is defined as $\mathbb{J}=N\times (N\cup\{\infty\})$, with the order defined by $(j,k)\leq (m,n)$ iff $j=m$ and $k\leq n$ or $n=\infty$ and $k\leq m$. Let $X=(\mathbb{J}, \tau_{\sigma})$. Then
\begin{enumerate}
\item[(\texttt{a})] $\mathbb{J}$ is a dcpo, thus $X$ is a $d$-space.
\item[(\texttt{b})] $Q(X)=Q_{s}(X)$.
\item[(\texttt{c})] $X$ is not well-filtered.
\item[(\texttt{d})] $X$ is not $s$-well-filtered.
\end{enumerate}
\end{example}

Using the equivalent classes of directed subsets, Ershov introduced one way to get $d$-completions of $T_{0}$ spaces in \cite{Ershov 17}. Inspired by his method, Liu, Li and Wu presented one way to get well-filtered reflections of $T_{0}$ spaces using the equivalent classes of KF-subsets in \cite{Liu2020}. Now for a $K$-subset system $Q_{k}:\mathbf{Top}_{0}\longrightarrow \mathbf{Set}$, we use the equivalent classes of $k$-Rudin sets to construct the $k$-well-filtered reflections of $T_{0}$ spaces. Let $(X,\tau)$ be a $T_{0}$ space. Consider an equivalence relation $\sim$ on $\mathrm{K}^{R}(X)$ which is defined as follows:
$$A_{0}\sim A_{1} \mbox{ if and only if } A_{0}\cap U\neq \emptyset \mbox{ is equivalent to } A_{1}\cap U\neq \emptyset \mbox{ for any } U\in \tau$$
where $A_{0},A_{1}\in \mathrm{K}^{R}(X)$. Note that $A_{0}\sim A_{1}$ if and only if $\mathrm{cl}_{X}(A_{0})=\mathrm{cl}_{X}(A_{1})$. Let
$$[A]=\{A^{\prime}\in \mathrm{K}^{R}(X)\mid A\sim A^{\prime}\}, A\in \mathrm{K}^{R}(X),$$
$$K(X)=\{[A]\mid A\in \mathrm{K}^{R}(X)\},$$
$$U^{\ast}=\{[A]\mid A\cap U\neq \emptyset\}, U\in \tau,$$
$$\tau^{\ast}=\{U^{\ast}\mid U\in \tau\}.$$

Then $\tau^{\ast}$ is a topology. Moreover, $(K(X), \tau^{\ast})$ is a $T_{0}$ space. For $Q_{k}(X)=S^{u}(X)$, we denote $K(X)$ by $D(X)$. For $Q_{k}(X)=Q(X)$, we denote $K(X)$ by $KF(X)$.

\begin{lemma}
Let $Q_{k}:\mathbf{Top}_{0}\longrightarrow \mathbf{Set}$ be a $K$-subset system, $(X,\tau)$ a $T_{0}$ space and $\lambda: X\longrightarrow K(X)$ the map defined by $\lambda (x)=[\{x\}]$. Then the map $\lambda$ is a homeomorphic embedding.
\end{lemma}

\begin{lemma}\label{4}
Let $Q_{k}:\mathbf{Top}_{0}\longrightarrow \mathbf{Set}$ be a $K$-subset system and $(X,\tau)$ a $T_{0}$ space. Then the following are equivalent:
\begin{enumerate}
\item[\texttt{1}.] X is a $k$-well-filtered space.
\item[\texttt{2}.] $K(X)\cong X$ (under the map $\lambda$).
\end{enumerate}
\end{lemma}

Let $Q_{k}:\mathbf{Top}_{0}\longrightarrow \mathbf{Set}$ be a $K$-subset system and $(X,\tau)$ a $T_{0}$ space. Suppose that $Y$ is a well-filtered space that has $X$ as a subspace. Since $Y$ is well-filtered, it is $k$-well-filtered. By~Lemma~\ref{4}, we have $K(Y)\cong Y$. In general, we can consider $K_{\gamma}(X)$ as a subspace of $K_{\beta}(X)$ in the sense of embedding mappings for all ordinals $\gamma\leq\beta$. The transfinite sequence of extensions is constructed as follows:
\begin{enumerate}
\item[(\texttt{1})] $K_{0}(X)=X$,
\item[(\texttt{2})] $K_{\beta+1}(X)=K(K_{\beta}(X))$,
\item[(\texttt{3})] $K_{\beta}(X)=\bigcup_{\gamma<\beta}K_{\gamma}(X)$  if $\beta$ is a limit ordinal.
\end{enumerate}

By \cite{Ershov 99} and \cite{Shen19}, we have the following similar results.

\begin{theorem}\label{$k$-well-filterification}
 For a $K$-subset system $Q_{k}:\mathbf{Top}_{0}\longrightarrow \mathbf{Set}$ and a $T_{0}$ space $(X,\tau)$, the $k$-well-filterification of $X$ exists; i.e., there exists an ordinal $\alpha$ such that $\mathrm{H}_{k}(X)=K_{\alpha}(X)\cong K_{\alpha+1}(X)$.
\end{theorem}
\begin{proof}
The proof is similar to the method of constructing the d-completion of $T_{0}$ spaces in \cite{Ershov 99} and the method of constructing the well-filterification of $T_{0}$ spaces in \cite{Shen19}.
\end{proof}

\begin{definition}
Let $Q_{k}:\mathbf{Top}_{0}\longrightarrow \mathbf{Set}$ be a $K$-subset system and $(X,\tau)$ a $T_{0}$ space. The \emph{$k$-rank} of $X$ is the least ordinal $\alpha$ such that $K_{\alpha}(X)\cong K_{\alpha+1}(X)$. We denote the $k$-rank of a space $X$ by $\mathrm{rank}_{k}(X)$.
\end{definition}

Similarly, the \emph{$d$-rank} of $X$ is the least ordinal $\alpha$ such that $D_{\alpha}(X)\cong D_{\alpha+1}(X)$, it is denoted by $\mathrm{rank}_{d}(X)$ in \cite{Ershov 17}. The \emph{$wf$-rank} of $X$ is the least ordinal $\alpha$ such that $KF_{\alpha}(X)\cong KF_{\alpha+1}(X)$, it is denoted by $\mathrm{rank}_{wf}(X)$ in \cite{Liu21}.

\section{$\alpha^{k}$-special spaces}
For a $K$-subset system $Q_{k}: \mathbf{Top}_{0}\longrightarrow \mathbf{Set}$, in Theorem \ref{$k$-well-filterification}, there exists an ordinal $\alpha$ such that $\mathrm{rank}_{k}(X)=\alpha$ for a $T_{0}$ space $X$. Conversely, for any given ordinal $\alpha$ it is natural to ask whether there exists a $T_{0}$ space $X$ such that $\mathrm{rank}_{k}(X)=\alpha$. In this section, we  prove that for any given ordinal $\alpha$, there exists a $T_{0}$ space $X$ such that $\mathrm{rank}_{k}(X)=\alpha$.

\begin{definition}\label{special}
Let $Q_{k}:\mathbf{Top}_{0}\longrightarrow \mathbf{Set}$ be a $K$-subset system and $(X,\tau)$ a $T_{0}$ space. For an ordinal $\alpha$, $X$ is called \emph{$\alpha^{k}$-special} if the following conditions are satisfied:
\begin{enumerate}
\item[(\texttt{1})] $\mathrm{rank}_{k}(X)=\alpha$;
\item[(\texttt{2})] $\alpha$ is the least ordinal for which $K_{\alpha}(X)$ has a greatest element.
\end{enumerate}
\end{definition}

$X$ is called \emph{$\alpha^{d}$-special} (resp., \emph{$\alpha^{wf}$-special}), similarly, see \cite{Ershov 17} and \cite{Liu21}, respectively.

\begin{remark}
If $X$ is a $\alpha^{k}$-special space, then $\alpha$ is not a limit ordinal.
\end{remark}

\begin{proof}
In fact, let $(X,\tau)$ be a $\alpha^{k}$-special space. Suppose that $\alpha$ is a limit ordinal. Then $K_{\alpha}(X)=\bigcup_{\beta<\alpha}K_{\beta}(X)$. By Definition \ref{special}, $K_{\alpha}(X)$ has a greatest element. Hence, there exists $\beta<\alpha$ such that $K_{\beta}(X)$ has a greatest element, which is a contradiction. So $\alpha$ will not be a limit ordinal.
\end{proof}

\begin{lemma}
For any nonlimit ordinal $\alpha$, every $\alpha^{k}$-special space is irreducible.
\end{lemma}

\begin{proof}
The proof is similar to Lemma 3.3 in \cite{Ershov 17}.
\end{proof}

Recall the following construction in \cite{Ershov 17}. For topological spaces $X$ and $Y_{x}$, $x\in X$, let
$$Z=\bigcup\limits_{x\in X}Y_{x}\times\{x\},$$
$$\tau=\{U\subseteq Z\mid (U)_{x}\in \tau (Y_{x}) \mbox{ for any }x\in X \mbox{ and }(U)_{X}\in \tau(X)\},$$
where $(U)_{x}=\{y\in Y_{x}\mid (y,x)\in U\}$ for any $x\in X$ and $(U)_{X}=\{x\in X \mid (U)_{x}\neq \emptyset\}$.

\begin{lemma}(\cite{Ershov 17})
Let $X$ be a $T_{0}$ space and $Y_{x}$ an irreducible $T_{0}$ space for every $x\in X$. Then
\begin{enumerate}
\item[(\texttt{1})] $\tau$ is a $T_{0}$ separable topology on $Z$.
\item[(\texttt{2})] The map $y\mapsto(y,x)$ determines a homeomorphic embedding of $Y_{x}$ in $Z$ for any $x\in X$.
\item[(\texttt{3})] If the space $X$ is irreducible, then the space $Z$ is also irreducible.
\end{enumerate}
\end{lemma}

For any subset $A\subseteq Z$, put $\widetilde{X}=\{x\in X\mid Y_{x} \mbox{ has the greatest element }\top_{x}\}$ with the induced topology of $X$. Define
$$(A)_{x}=\{y\in Y_{x}\mid (y,x)\in A\} \mbox{ for any }x\in X,$$
$$(A)_{X}=\{x\in X\mid (A)_{x}\neq \emptyset\},$$
$$A_{\ast}=\{x\in \widetilde{X}\mid (\top_{x},x)\in A\}.$$

And the space $(Z, \tau)$ is also denoted by $\sum\limits_{X}Y_{x}$.

\begin{lemma}(\cite{Ershov 17})\label{order}
Let $X$ be a $T_{0}$ space, $Y_{x}$ an irreducible $T_{0}$ space for any $x\in X$, and $Z=\sum\limits_{X}Y_{x}$. For all $(y_{0},x_{0})$, $(y_{1},x_{1})\in Z$, we have $(y_{0},x_{0})\leq(y_{1},x_{1})$ if and only if the following two alternatives hold:
\begin{enumerate}
\item[(\texttt{1})] $x_{0}=x_{1}$ and $y_{0}\leq_{Y_{x_{0}}}y_{1}$;
\item[(\texttt{2})] $x_{0}<_{X}x_{1}$ and $y_{1}=\top_{x_{1}}$ is the greatest element in $Y_{x_{1}}$.
\end{enumerate}
\end{lemma}

\begin{lemma}(\cite{Ershov 17})
Let $X$ be a $T_{0}$ space, $Y_{x}$ an irreducible $T_{0}$ space for any $x\in X$, and $Z=\sum\limits_{X}Y_{x}$. Then an arbitrary set $S^{\prime}\in \mathcal{D}(Z)$ contains a cofinal subset $S\subseteq S^{\prime}$ having one of the following forms:
\begin{enumerate}
\item[(\texttt{\expandafter\romannumeral01})] $S=\{(y,x)\mid y\in (S)_{x}\}$ for some fixed $x\in X$ and $(S)_{x}\in \mathcal{D}(Y_{x})$;
\item[(\texttt{\expandafter\romannumeral02})] $S=\{(\top_{x},x)\mid x\in S_{\ast}\}$ for some $S_{\ast}\in \mathcal{D}(\widetilde{X})$.
\end{enumerate}
\end{lemma}

For any irreducible topological space $Y$, put
\[Y^{\top}=
 \begin{cases}
 Y,& \mbox{ if } Y \mbox{ has a greatest element},\\

\langle Y\cup\{\top\}, \tau(Y)^{\top} \rangle,& \mbox{ otherwise},\\

\end{cases}\]
where $\tau(Y)^{\top}=\{U\cup \{\top\}\mid \emptyset \neq U\in \tau(Y)\}\cup\{\varnothing\}$. It is easy to see that for any irreducible $T_{0}$ space $Y$, $Y^{\top}$ is also a $T_{0}$ space and has a greatest element. Let
$$X^{\prime}=\{[\{x\}]\mid x\in X\}\cup D(\widetilde{X})\subseteq D(X).$$
Then $X^{\prime}$ with the induced topology is a subspace of $D(X)$ and the space $D(Y_{x})$ is irreducible for any $x\in X$ from \cite{Ershov 17}. Moreover, for any $x^{\prime}\in X^{\prime}$, we define
\[Y^{\prime}_{x^{\prime}}=
 \begin{cases}
D(Y_{x}),& \mbox{ if } x^{\prime}=[\{x\}] \mbox{ for some $x\in X$},\\

\top,& \mbox{ otherwise},\\

\end{cases}\]
where $\top=\langle\{\top\},\{\emptyset,\{\top\}\}\rangle$. Then we have the following theorem.

\begin{theorem}(\cite{Ershov 17})\label{X and Xpie}
Let $X$ be a $T_{0}$ space, $Y_{x}$ an irreducible $T_{0}$ space for any $x\in X$, and $Z=\sum\limits_{X}Y_{x}$. Then the spaces $D(Z)$ and $Z^{\prime}=\sum\limits_{X^{\prime}}Y^{\prime}_{X^{\prime}}$ are homeomorphic.
\end{theorem}

For any ordinal $\alpha>0$, consider the irreducible $T_{0}$ space
$$\mathbb{O}_{\alpha}=\langle\da \alpha\setminus \{\alpha\}, \{\emptyset\}\cup\{\ua\beta \mid \beta<\alpha \mbox{ is not a limit ordinal}\}\rangle.$$

\begin{proposition}(\cite{Ershov 17})\label{Ershov 17}
Let $\alpha> 0$ be an ordinal.
\begin{enumerate}
\item[(\texttt{\expandafter\romannumeral01})] If $\alpha$ is a limit ordinal, then $\mathrm{H}_{d}(\mathbb{O}_{\alpha})=\mathbb{O}^{\top}_{\alpha}=D(\mathbb{O}_{\alpha})$, i.e., the space $\mathbb{O}_{\alpha}$ is $1^{d}$-special.
\item[(\texttt{\expandafter\romannumeral02})] If $\alpha$ is not a limit ordinal, then $\mathrm{H}_{d}(\mathbb{O}_{\alpha})=\mathbb{O}_{\alpha}$, i.e., the $d$-rank of $\mathbb{O}_{\alpha}$ is equal to $0$.
\item[(\texttt{\expandafter\romannumeral03})] If $\alpha$ is a limit ordinal, $\gamma$ is not a limit ordinal, a $T_{0}$ space $Y_{\beta}$ is $\gamma^{d}$-special for any $\beta< \alpha$ and the space $Z=\sum\limits_{\mathbb{O}_{\alpha}}Y_{\beta}$, then the spaces  $D_{\delta}(Z)\cong\sum\limits_{\mathbb{O}_{\alpha}}D_{\delta}(Y_{\beta})$ for any ordinal $\delta\leq \gamma$.
\item[(\texttt{\expandafter\romannumeral04})] If $\alpha$ is a limit ordinal, $\gamma$ is not a limit ordinal and a $T_{0}$ space $Y_{\beta}$ is $\gamma^{d}$-special for any $\beta< \alpha$, then the space $Z=\sum\limits_{\mathbb{O}_{\alpha}}Y_{\beta}$ is $(\gamma+1)^{d}$-special.
\item[(\texttt{\expandafter\romannumeral05})] If $Y$ is an $(\alpha+1)^{d}$-special for some ordinal $\alpha$, then $D_{\beta}(Y^{\top})\cong D_{\beta}(Y)^{\top}$ for any $\beta\leq\alpha$ and $D_{\alpha+1}(Y^{\top})\cong D_{\alpha+1}(Y)=\mathrm{H}_{d}(Y)$.
\item[(\texttt{\expandafter\romannumeral06})] If $\alpha$ is a limit ordinal, and a $T_{0}$ space $Y_{\beta}$ is $(\beta+1)^{d}$-special for any $\beta< \alpha$, then the space $Z=\sum\limits_{\mathbb{O}_{\alpha}}Y_{\beta}$ is $(\alpha+1)^{d}$-special and the $d$-rank of a space $Z^{\top}$ is equal to $\alpha$.
\end{enumerate}
\end{proposition}

\begin{lemma}
Let $\mathbb{N}=(N, \tau_{\sigma})$ denote the set $N$ of natural numbers endowed with the Scott topology.  Then $\mathbb{N}$ is $1^{d}$-special.
\end{lemma}

\begin{proof}
This directly follows from ($\expandafter{\mathrm{\romannumeral01}}$) of Proposition \ref{Ershov 17}.
\end{proof}

\begin{lemma}\label{d rank}
Let $\alpha> 0$ be an ordinal.
\begin{enumerate}
\item[(\texttt{1})] If $\alpha$ is not a limit ordinal and a $T_{0}$ space $X_{n}$ is $\alpha^{d}$-special for any $n\in N$, then the space $Z=\sum\limits_{\mathbb{N}}X_{n}$ is $(\alpha+1)^{d}$-special.
\item[(\texttt{2})] If $\alpha$ is a limit ordinal and a $T_{0}$ space $X_{n}$ is $(\overline{\alpha}+n)^{d}$-special for any $n\in N$, then the space $Z=\sum\limits_{\mathbb{N}}X_{n}$ is $(\alpha+1)^{d}$-special and the $d$-rank of the space $Z^{\top}$ is equal to $\alpha$, where $\overline{\alpha}=0$ if $\alpha= \omega$, otherwise, $\overline{\alpha}$ denotes the largest limit ordinal less than $\alpha$.
\end{enumerate}
\end{lemma}

\begin{proof}
\begin{enumerate}
\item[(\texttt{1})] It follows directly from ($\expandafter{\mathrm{\romannumeral04}}$) of Proposition \ref{Ershov 17}.
\item[(\texttt{2})] First we prove that the spaces $D_{\delta}(Z)$ and $\sum\limits_{\mathbb{N}}W^{\delta}_{n}$ are homeomorphic for every ordinal $\delta<\alpha$, where
    \[W^{\delta}_{n}=
 \begin{cases}
 D_{\delta}(X_{n}),& \mbox{if } \delta<\overline{\alpha}+n,\\

 \mathrm{H}_{d}(X_{n}),& \mbox{if } \overline{\alpha}+n\leq\delta<\alpha. \\

\end{cases}\]
We use induction on $\delta$. For $\delta=0$, the statement follows from the definition of space $Z$.

Let $\delta$ be an ordinal such that $\delta+1<\alpha$, and suppose that $D_{\delta}(Z)\cong\sum\limits_{\mathbb{N}}W^{\delta}_{n}$. Then the space $D_{\delta}(X_{n})$ does not contain a greatest element for the arbitrary $n\in N$ such that $\delta<\overline{\alpha}+n$. Note that $\widetilde{N}=\{n\in N \mid \overline{\alpha}+n\leq \delta\}$ is a finite subset in $N$, hence $\mathbb{N}^{\prime}=\{[\{n\}]\mid n\in N\}\cup D(\widetilde{N})\cong\mathbb{N}$. By Theorem \ref{X and Xpie}, we have
    $$D_{\delta+1}(Z)=D(D_{\delta}(Z))\cong D(\sum\limits_{\mathbb{N}}W^{\delta}_{n})\cong \sum\limits_{\mathbb{N}}D(W^{\delta}_{n})=
    \sum\limits_{\mathbb{N}}W^{\delta+1}_{n}.$$
Suppose now that $\delta<\alpha$ is a limit ordinal and $D_{\beta}(Z)\cong\sum\limits_{\mathbb{N}}W^{\beta}_{n}$ for any ordinal $\beta<\delta$. By~Theorem \ref{X and Xpie}, we get
$$D_{\delta}(Z)=\bigcup\limits_{\beta<\delta}D_{\beta}(Z)\cong \bigcup\limits_{\beta<\delta}\sum\limits_{\mathbb{N}}W^{\beta}_{n}\cong \sum\limits_{\mathbb{N}}\bigcup\limits_{\beta<\delta}W^{\beta}_{n}=
\sum\limits_{\mathbb{N}}W^{\delta}_{n}.$$
Thus by induction, we have $D_{\delta}(Z)\cong \sum\limits_{\mathbb{N}}W^{\delta}_{n}$ for any ordinal $\delta<\alpha$.

Therefore,
$$D_{\alpha}(Z)=\bigcup\limits_{\delta<\alpha}D_{\delta}(Z)\cong \bigcup\limits_{\delta<\alpha}\sum\limits_{\mathbb{N}}W^{\delta}_{n}\cong \sum\limits_{\mathbb{N}}\bigcup\limits_{\delta<\alpha}W^{\delta}_{n}=
\sum\limits_{\mathbb{N}}W^{\alpha}_{n}\cong
\sum\limits_{\mathbb{N}}\mathrm{H}_{d}(X_{n}).$$
Moreover, for any $n\in N$, the space $\mathrm{H}_{d}(X_{n})$ has a greatest element, which implies that $\widetilde{N}=N$. Hence $\mathbb{N}^{\prime}\cong\mathbb{N}^{\top}$. In the view of Lemma \ref{order} and Theorem \ref{X and Xpie}, we obtain
$$D_{\alpha+1}(Z)=D(D_{\alpha}(Z))\cong D(\sum\limits_{\mathbb{N}}\mathrm{H}_{d}(X_{n}))\cong \sum\limits_{\mathbb{N}^{\top}}X^{\prime}_{n^{\prime}}\cong
(\sum\limits_{\mathbb{N}}\mathrm{H}_{d}(X_{n}))^{\top},$$
where \[X^{\prime}_{n^{\prime}}=
 \begin{cases}

 \mathrm{H}_{d}(X_{n^{\prime}}),& \mbox{if }n^{\prime}\in N,\\

 \mathbb{\top},& \mbox{if }n^{\prime}=\top,\\

 \end{cases}\]
and
$$D_{\alpha+2}(Z)=D(D_{\alpha+1}(Z))\cong D(\sum\limits_{\mathbb{N}^{\top}}X^{\prime}_{n^{\prime}})\cong \sum\limits_{\mathbb{N}^{\top}}D(X^{\prime}_{n^{\prime}})\cong
\sum\limits_{\mathbb{N}^{\top}}X^{\prime}_{n^{\prime}}\cong D_{\alpha+1}(Z).$$

Again by Lemma \ref{order}, the space $D_{\beta}(Z)$ has not a greatest element for any ordinal $\beta\leq\alpha$. Therefore, by virtue of Definition \ref{special}, the space $Z$ is $(\alpha+1)^{d}$-special.

For $Z^{\top}$, first, we claim that the spaces $D_{\delta}(Z^{\top})\cong\sum\limits_{\mathbb{N^{\top}}}W^{\delta}_{n^{\prime}}$ for the arbitrary ordinal $\delta\leq\alpha$, where
    \[W^{\delta}_{n^{\prime}}=
 \begin{cases}
 D_{\delta}(X_{n^{\prime}}),& \mbox{if } \delta<\overline{\alpha}+n^{\prime}<\alpha,\\

 \mathrm{H}_{d}(X_{n^{\prime}}),& \mbox{if } \overline{\alpha}+n^{\prime}\leq\delta\leq\alpha, \\

 \mathbb{\top},& \mbox{if } n^{\prime}=\top, \\

\end{cases}\]

By the part $($\expandafter{\romannumeral05}$)$ of Proposition \ref{Ershov 17}, we get
$$D_{\delta}(Z^{\top})\cong (D_{\delta}(Z))^{\top}\cong (\sum\limits_{\mathbb{N}}W^{\delta}_{n})^{\top}\cong \sum\limits_{\mathbb{N^{\top}}}W^{\delta}_{n^{\prime}}$$
for every ordinal $\delta\leq\alpha$. This implies that $$D_{\alpha}(Z^{\top})\cong \sum\limits_{\mathbb{N^{\top}}}W^{\alpha}_{n^{\prime}}\cong
(\sum\limits_{\mathbb{N}}\mathrm{H}_{d}(X_{n}))^{\top}\cong D_{\alpha+1}(Z),$$
which is a $d$-space by the above proof. Therefore, $D_{\alpha}(Z^{\top})\cong D_{\alpha+1}(Z^{\top})$.

Next, we claim that $D_{\delta}(Z^{\top})$ is not a $d$-space for any ordinal $\delta<\alpha$. Assume that there exists an ordinal $\delta<\alpha$ such that $D_{\delta}(Z^{\top})$ is a $d$-space. Then by Lemma \ref{4}, $D_{\delta}(Z^{\top})\cong D_{\alpha}(Z^{\top})$. However, from the preceding discussion, we have that the spaces $D_{\delta}(Z^{\top})\cong\sum\limits_{\mathbb{N^{\top}}}W^{\delta}_{n^{\prime}}$ for the arbitrary ordinal $\delta<\alpha$. Note that there are at most finitely many $W^{\delta}_{n^{\prime}}$'s are $d$-spaces. Furthermore, for $\delta<\overline{\alpha}+n^{\prime}<\alpha$, $W^{\delta}_{n^{\prime}}=D_{\delta}(X_{n^{\prime}})$ is not a $d$-space and $W^{\alpha}_{n^{\prime}}=\mathrm{H}_{d}(X_{n^{\prime}})$ is a $d$-space, which implies that $W^{\delta}_{n^{\prime}}$ and $W^{\alpha}_{n^{\prime}}$ are not homeomorphic. Hence, $D_{\delta}(Z^{\top})$ and $D_{\alpha}(Z^{\top})$ are not homeomorphic, which is a contradiction. So the $d$-rank of the space $Z^{\top}$ is equal to $\alpha$.
\end{enumerate}
\end{proof}
For the $wf$-rank of a $T_{0}$ space, we have the following similar results.

\begin{lemma}\label{Liu bei}\cite{Liu22}
Let $X$ be a $T_{0}$ space, $Y_{x}$ an irreducible $T_{0}$ space for any $x\in X$, and $Z=\sum\limits_{X}Y_{x}$. Then an arbitrary set $A^{\prime}\in \overline{KF}(Z) $ contains a subset $A\subseteq A^{\prime}$ such that $A\sim A^{\prime}$ having one of the following forms:
\begin{enumerate}
\item[ (1)] 
there exists an element $x\in X$ such that $A\subseteq Y_{x}\times\{x\}$ and $(A)_{x}\in \overline{KF}(Y_{x})$;
\item[ (2)]
$A=\{(\top_{x},x)\mid x\in A_{\ast}\}$ for some $A_{\ast}\in \overline{KF}(\widetilde{X})$, where $\widetilde{X}=\{x\in X\mid Y_{x} \mbox{ has a greatest element}\}$.
\end{enumerate}
\end{lemma}

\begin{theorem}\label{WF X and Xpie}
Let $X$ be a $T_{0}$ space, $Y_{x}$ an irreducible $T_{0}$ space for any $x\in X$, and $Z=\sum\limits_{X}Y_{x}$. Then the spaces $KF(Z)$ and $Z^{\prime}=\sum\limits_{X^{\prime}}Y^{\prime}_{X^{\prime}}$ are homeomorphic, where
$$X^{\prime}=\{[\{x\}]\mid x\in X\}\cup KF(\widetilde{X})\subseteq KF(X),$$
and for any $x^{\prime}\in X^{\prime}$
\[Y^{\prime}_{x^{\prime}}=
 \begin{cases}
KF(Y_{x}),& \mbox{ if } x^{\prime}=[\{x\}] \mbox{ for some $x\in X$},\\

\top,& \mbox{ otherwise}, \mbox{ where } \top=\langle\{\top\},\{\emptyset,\{\top\}\}\rangle.\\

\end{cases}\]
\end{theorem}

\begin{lemma}\label{wf rank 1}
$(N, \tau_{\sigma})$ is $1^{wf}$-special.
\end{lemma}
\begin{lemma}\label{wf rank 2}
Let $\alpha> 0$ be an ordinal.
\begin{enumerate}
\item[(\texttt{1})] If $\alpha$ is not a limit ordinal and a $T_{0}$ space $X_{n}$ is $\alpha^{wf}$-special for any $n\in N$, then the space $Z=\sum\limits_{\mathbb{N}}X_{n}$ is $(\alpha+1)^{wf}$-special.
\item[(\texttt{2})] If $\alpha$ is a limit ordinal and a $T_{0}$ space $X_{n}$ is $(\overline{\alpha}+n)^{wf}$-special for any $n\in N$, then the space $Z=\sum\limits_{\mathbb{N}}X_{n}$ is $(\alpha+1)^{wf}$-special and the $wf$-rank of the space $Z^{\top}$ is equal to $\alpha$, where $\overline{\alpha}=0$ if $\alpha= \omega$, otherwise, $\overline{\alpha}$ denotes the largest limit ordinal less than $\alpha$.
\end{enumerate}
\end{lemma}

Next, let $Q_{k}$ be a $K$-subset system and $X$ a $T_{0}$ space. We denote

$$\overline{D}(X)=\{A\subseteq X\mid \mbox{ there exists a directed subset }D \mbox{ in }X \mbox{ such that }\overline{A}=\overline{D}\}.$$

For the $k$-rank of $X$, we deduce the following results.

\begin{lemma}\label{N is 1}
For $\mathbb{N}=(N, \tau_{\sigma})$, the following statements hold:
\begin{enumerate}
\item[(\texttt{1})] $\mathbb{N}$ is $1^{d}$-special and $1^{wf}$-special.
\item[(\texttt{2})] $\overline{D}(\mathbb{N})=\overline{KF}(\mathbb{N})$ and $D(\mathbb{N})=K(\mathbb{N})=KF(\mathbb{N})$.
\item[(\texttt{3})] $\mathbb{N}$ is $1^{k}$-special.
\end{enumerate}
\end{lemma}

\begin{lemma}\label{not a limit ordinal}
If $\alpha$ is not a limit ordinal and a $T_{0}$ space $X_{n}$ satisfies the following conditions.
    \begin{enumerate}
\item[(\texttt{1})] $X_{n}$ is $\alpha^{d}$-special and $\alpha^{wf}$-special,
\item[(\texttt{2})] $\overline{D}(K_{m}(X_{n}))=\overline{KF}(K_{m}(X_{n}))$ and $D_{m+1}(X_{n})=K_{m+1}(X_{n})=KF_{m+1}(X_{n})$ for $0\leq m<\alpha$,
\end{enumerate}
for any $n\in N$, then the space $Z=\sum\limits_{\mathbb{N}}X_{n}$ satisfies:
\begin{enumerate}
\item[(\texttt{1})] $Z$ is $(\alpha+1)^{d}$-special and $(\alpha+1)^{wf}$-special.
\item[(\texttt{2})] $\overline{D}(K_{m}(Z))=\overline{KF}(K_{m}(Z))$ and $D_{m+1}(Z)=K_{m+1}(Z)=KF_{m+1}(Z)$ for $0\leq m<\alpha+1$.
\item[(\texttt{3})] $Z$ is $(\alpha+1)^{k}$-special.
\end{enumerate}
\end{lemma}

\begin{proof}
For (1), this directly follows from Lemma \ref{d rank} (1) and Lemma \ref{wf rank 2} (1).

For (2), the proof is by induction on $m$.

Basic steps. For $m=0$, obviously, $\overline{D}(Z)\subseteq \overline{KF}(Z)$. Conversely, let $A\in \overline{KF}(Z)$. From Lemma~\ref{Liu bei}, there exists a subset $A^{\prime}\subseteq A$ that $A^{\prime}\sim A$ and $A^{\prime}$ satisfies Type (\expandafter{\romannumeral01}) in Lemma \ref{Liu bei}. This means that there exists $n\in N$ such that $A^{\prime}\subseteq X_{n}\times\{n\}$ and $(A^{\prime})_{n}\in \overline{KF}(X_{n})=\overline{D}(X_{n})$. By the definition of $\overline{D}(X_{n})$, there is a directed subset $D$ in $X_{n}$ such that $\mathrm{cl}_{X_{n}}((A^{\prime})_{n})=\mathrm{cl}_{X_{n}}(D)$. We claim that $\mathrm{cl}_{Z}(A^{\prime})=\mathrm{cl}_{Z}(D\times\{n\})$. For $(y,x)\in \mathrm{cl}_{Z}(A^{\prime})$, let $U$ be an open neighbourhood of $(y,x)$. Then we have that $U\cap A^{\prime}\neq \emptyset$. This implies that $(U)_{n}\cap (A^{\prime})_{n}\neq \emptyset$. Since $U$ is open in $Z$, we have $(U)_{n}\in \tau(X_{n})$. By $\mathrm{cl}_{X_{n}}((A^{\prime})_{n})=\mathrm{cl}_{X_{n}}(D)$, we have $(U)_{n}\cap D\neq \emptyset$, that is $U\cap (D\times \{n\})\neq \emptyset$. Hence $(y,x)\in \mathrm{cl}_{Z}(D\times\{n\})$. The opposite direction is similar to prove. So $A\in \overline{D}(Z)$. That is $\overline{D}(Z)= \overline{KF}(Z)$, which implies that $\overline{D}(Z)=\overline{K}(Z)= \overline{KF}(Z)$. Therefore, $D(Z)=K(Z)=KF(Z)$.

Inductive steps. There are two cases to consider:

Case 1. Let $m$ be an ordinal such that $m+1<\alpha+1$. Assume that $\overline{D}(K_{m}(Z))=\overline{KF}(K_{m}(Z))$ and $D_{m+1}(Z)=K_{m+1}(Z)=KF_{m+1}(Z)$, by Lemma \ref{Ershov 17} (\expandafter\romannumeral03), we have

$$D_{m+1}(Z)=K_{m+1}(Z)=KF_{m+1}(Z)\cong \sum\limits_{\mathbb{N}}D_{m+1}(X_{n}).$$
To prove $\overline{D}(K_{m+1}(Z))=\overline{KF}(K_{m+1}(Z))$, it is enough to show $\overline{D}(\sum\limits_{\mathbb{N}}D_{m+1}(X_{n}))=
\overline{KF}(\sum\limits_{\mathbb{N}}D_{m+1}(X_{n}))$. Clearly, $\overline{D}(\sum\limits_{\mathbb{N}}D_{m+1}(X_{n}))\subseteq
\overline{KF}(\sum\limits_{\mathbb{N}}D_{m+1}(X_{n}))$. Conversely, for any $A\in \overline{KF}(\sum\limits_{\mathbb{N}}D_{m+1}(X_{n}))$, by Lemma~\ref{Liu bei}, there exists a subset $A^{\prime}\subseteq A$ such that $A^{\prime}\sim A$. Two options are possible:

Case 1.1. $A^{\prime}$ is Type (\expandafter\romannumeral01) in Lemma~\ref{Liu bei}. This implies that there exists $n\in N$ such that $A^{\prime}\subseteq D_{m+1}(X_{n})\times \{n\}$ and $(A^{\prime})_{n}\in \overline{KF}(D_{m+1}(X_{n}))$. By the condition (2) of $X_{n}$, we get

$$D_{m+1}(X_{n})=K_{m+1}(X_{n})=KF_{m+1}(X_{n})\mbox{ and } \overline{KF}(K_{m+1}(X_{n}))=\overline{D}(K_{m+1}(X_{n})).$$
Hence, $(A^{\prime})_{n}\in \overline{KF}(D_{m+1}(X_{n}))=\overline{D}(D_{m+1}(X_{n}))$. So $A^{\prime}\in \overline{D}(\sum\limits_{\mathbb{N}}(D_{m+1}(X_{n})))$. This implies that $A\in \overline{D}(\sum\limits_{\mathbb{N}}(D_{m+1}(X_{n})))$.

Case 1.2. $A^{\prime}$ is Type (\expandafter\romannumeral02) in Lemma~\ref{Liu bei}. This means that
there exists $A_{\ast}\in \overline{KF}(\widetilde{\mathbb{N}})$ such that $A^{\prime}=\{(\top_{n},n)\mid n\in A_{\ast}\}$. Note that $m+1=\alpha$, then $\sum\limits_{\mathbb{N}}D_{m+1}(X_{n})\cong
\sum\limits_{\mathbb{N}}\mathrm{H}_{d}(X_{n})\cong
\sum\limits_{\mathbb{N}}\mathrm{H}_{wf}(X_{n})$. So $\widetilde{\mathbb{N}}=\mathbb{N}$. Hence, $\overline{KF}(\widetilde{\mathbb{N}})=\overline{D}(\widetilde{\mathbb{N}})$, which implies that $A^{\prime}\in \overline{D}(\sum\limits_{\mathbb{N}}D_{m+1}(X_{n}))$. Therefore, $A\in \overline{D}(\sum\limits_{\mathbb{N}}D_{m+1}(X_{n}))$.

In any case, we have $\overline{KF}(\sum\limits_{\mathbb{N}}D_{m+1}(X_{n}))\subseteq
\overline{D}(\sum\limits_{\mathbb{N}}D_{m+1}(X_{n}))
$. So $\overline{D}(K_{m+1}(Z))=\overline{KF}(K_{m+1}(Z))$. This implies that $D_{m+2}(Z)=K_{m+2}(Z)=KF_{m+2}(Z)$.

Case 2. Suppose that $m<\alpha+1$ is a limit ordinal and the required statement holds for any $\delta<m$. Then
$$K_{m}(Z)=\bigcup\limits_{\delta<m}K_{\delta}(Z)\cong
\bigcup\limits_{\delta<m}D_{\delta}(Z)\cong
\bigcup\limits_{\delta<m}\sum\limits_{\mathbb{N}}D_{\delta}(X_{n})\cong
\sum\limits_{\mathbb{N}}\bigcup\limits_{\delta<m}D_{\delta}(X_{n})\cong
\sum\limits_{\mathbb{N}}D_{m}(X_{n}).$$

Now it is enough to show that $\overline{KF}(\sum\limits_{\mathbb{N}}D_{m}(X_{n}))=
\overline{D}(\sum\limits_{\mathbb{N}}D_{m}(X_{n}))$. Repeat the proof method of Case 1, we get that $\overline{D}(K_{m}(Z))=\overline{KF}(K_{m}(Z))$. So $\overline{D}(K_{m}(Z))=\overline{K}(K_{m}(Z))=\overline{KF}(K_{m}(Z))$. Hence, $D_{m+1}(Z)=K_{m+1}(Z)=KF_{m+1}(Z)$.

For (3), by (2), let $m=\alpha$. We get $$D_{\alpha+1}(Z)=K_{\alpha+1}(Z)=KF_{\alpha+1}(Z).$$
Since $Z$ is $(\alpha+1)^{d}$-special and $(\alpha+1)^{wf}$-special, for any ordinal $\delta<\alpha+1$, $K_{\delta}(Z)$ is not a $d$-space and $K_{\alpha+1}(Z)$ is well-filtered. Therefore, $K_{\delta}(Z)$ is not $k$-well-filtered and $K_{\alpha+1}(Z)$ is $k$-well-filtered. Then $Z$ is $(\alpha+1)^{k}$-special.
\end{proof}

\begin{lemma}\label{limit rank}
If $\alpha$ is a limit ordinal and a $T_{0}$ space $X_{n}$ satisfies the following conditions.
    \begin{enumerate}
\item[(\texttt{1})] $X_{n}$ is $(\overline{\alpha}+n)^{d}$-special and $(\overline{\alpha}+n)^{wf}$-special,
\item[(\texttt{2})] $\overline{D}(K_{m}(X_{n}))=\overline{KF}(K_{m}(X_{n}))$ and $D_{m+1}(X_{n})=K_{m+1}(X_{n})=KF_{m+1}(X_{n})$ for $0\leq m<\overline{\alpha}+n$,
\end{enumerate}
for any $n\in N$, then for the space $Z=\sum\limits_{\mathbb{N}}X_{n}$, we have the following conclusions.
\begin{enumerate}
\item[(\texttt{1})] $Z$ is $(\alpha+1)^{d}$-special and $(\alpha+1)^{wf}$-special.
\item[(\texttt{2})] $\overline{D}(K_{m}(Z))=\overline{KF}(K_{m}(Z))$ and $D_{m+1}(Z)=K_{m+1}(Z)=KF_{m+1}(Z)$ for $0\leq m<\alpha+1$.
\item[(\texttt{3})] $Z$ is $(\alpha+1)^{k}$-special.
\end{enumerate}
Moreover, for the space $Z^{\top}$, the following results hold.
\begin{enumerate}
\item[(\texttt{1})] $rank_{d}(Z^{\top})=rank_{wf}(Z^{\top})=\alpha$,
\item[(\texttt{2})] $\overline{D}(K_{m}(Z^{\top}))=\overline{KF}(K_{m}(Z^{\top}))$ and $D_{m+1}(Z^{\top})=K_{m+1}(Z^{\top})=KF_{m+1}(Z^{\top})$ for $0\leq m<\alpha$,
\item[(\texttt{3})] $rank_{k}(Z^{\top})=\alpha$.
\end{enumerate}
\end{lemma}

\begin{proof}
First, we consider the space $Z=\sum\limits_{\mathbb{N}}X_{n}$.

For (1), it follows directly from Lemma \ref{d rank} (2) and Lemma \ref{wf rank 2} (2).

For (2), we proceed by induction. For $m=0$, the statement follows from the proof of Lemma \ref{not a limit ordinal}. Let $m$ be an ordinal such that $m+1<\alpha+1$. Assume that $\overline{D}(K_{m}(Z))=\overline{KF}(K_{m}(Z))$ and $D_{m+1}(Z)=K_{m+1}(Z)=KF_{m+1}(Z)$, by the proof of Lemma \ref{d rank} (2), we derive
$$D_{m+1}(Z)=K_{m+1}(Z)=KF_{m+1}(Z)\cong \sum\limits_{\mathbb{N}}W^{m+1}_{n}$$ where
    \[W^{m+1}_{n}=
 \begin{cases}
 D_{m+1}(X_{n}),& \mbox{if } m+1<\overline{\alpha}+n,\\

 \mathrm{H}_{d}(X_{n}),& \mbox{if } \overline{\alpha}+n\leq m+1<\alpha+1. \\

\end{cases}\]
To prove $\overline{D}(K_{m+1}(Z))=\overline{KF}(K_{m+1}(Z))$, it suffices to show that $\overline{D}(\sum\limits_{\mathbb{N}}W^{m+1}_{n})=
\overline{KF}(\sum\limits_{\mathbb{N}}W^{m+1}_{n})$. Clearly, $\overline{D}(\sum\limits_{\mathbb{N}}W^{m+1}_{n})\subseteq
\overline{KF}(\sum\limits_{\mathbb{N}}W^{m+1}_{n})$. For any $A\in \overline{KF}(\sum\limits_{\mathbb{N}}W^{m+1}_{n})$, by Lemma \ref{Liu bei}, there exists a subset $A^{\prime}\subseteq A$ such that $A^{\prime}\sim A$. There are two cases to consider:

Case 1. $A^{\prime}$ is Type (\expandafter\romannumeral01) in Lemma \ref{Liu bei}. This implies that there exists $n\in N$ such that $A^{\prime}\subseteq W^{m+1}_{n}\times \{n\}$ and $(A^{\prime})_{n}\in \overline{KF}(W^{m+1}_{n})$. Again there are two cases to consider:

Case 1.1. If $m+1<\overline{\alpha}+n$, then $W^{m+1}_{n}=D_{m+1}(X_{n})$. By the condition (2) of $X_{n}$, we have
$$D_{m+1}(X_{n})=K_{m+1}(X_{n})=KF_{m+1}(X_{n})\mbox{ and } \overline{KF}(K_{m+1}(X_{n}))=\overline{D}(K_{m+1}(X_{n})).$$
Hence, $(A^{\prime})_{n}\in \overline{KF}(D_{m+1}(X_{n}))=\overline{D}(D_{m+1}(X_{n}))$. So $A^{\prime}\in \overline{D}(\sum\limits_{\mathbb{N}}W^{m+1}_{n})$. This implies that $A\in \overline{D}(\sum\limits_{\mathbb{N}}W^{m+1}_{n})$.

Case 1.2. If $\overline{\alpha}+n\leq m+1<\alpha+1$, then $W^{m+1}_{n}=\mathrm{H}_{wf}(X_{n})=\mathrm{H}_{d}(X_{n})$. Hence, $\overline{KF}(W^{m+1}_{n})=\overline{D}(W^{m+1}_{n})$. It is straightforward to check that $A\in \overline{D}(\sum\limits_{\mathbb{N}}W^{m+1}_{n})$.

Case 2. $A^{\prime}$ is Type (\expandafter\romannumeral02) in Lemma \ref{Liu bei}. This means that
there exists $A_{\ast}\in \overline{KF}(\widetilde{\mathbb{N}})$ such that $A^{\prime}=\{(\top_{n},n)\mid n\in A_{\ast}\}$. Note that there are at most finitely many $W^{m+1}_{n}$'s which have a greatest element; thus $\overline{KF}(\widetilde{\mathbb{N}})=\overline{D}(\widetilde{\mathbb{N}})$. This implies that $A^{\prime}\in \overline{D}(\sum\limits_{\mathbb{N}}W^{m+1}_{n})$. Therefore, $A\in \overline{D}(\sum\limits_{\mathbb{N}}W^{m+1}_{n})$.

In any case, we have that $\overline{KF}(\sum\limits_{\mathbb{N}}W^{m+1}_{n})\subseteq
\overline{D}(\sum\limits_{\mathbb{N}}W^{m+1}_{n})
$. So $\overline{D}(K_{m+1}(Z))=\overline{KF}(K_{m+1}(Z))$. This implies that $D_{m+2}(Z)=K_{m+2}(Z)=KF_{m+2}(Z)$.

Suppose that $m<\alpha+1$ is a limit ordinal and that the required statement holds for any $\delta<m$. Then
$$K_{m}(Z)=\bigcup\limits_{\delta<m}K_{\delta}(Z)\cong
\bigcup\limits_{\delta<m}D_{\delta}(Z)\cong
\bigcup\limits_{\delta<m}\sum\limits_{\mathbb{N}}W^{\delta}_{n}\cong
\sum\limits_{\mathbb{N}}\bigcup\limits_{\delta<m}W^{\delta}_{n}\cong
\sum\limits_{\mathbb{N}}W^{m}_{n}$$

where \[W^{\delta}_{n}=
 \begin{cases}
 D_{\delta}(X_{n}),& \mbox{if } \delta<\overline{\alpha}+n,\\

 \mathrm{H}_{d}(X_{n}),& \mbox{if } \overline{\alpha}+n\leq\delta<\alpha. \\

\end{cases}\]
Now it is enough to show that $\overline{KF}(\sum\limits_{\mathbb{N}}W^{m}_{n})=
\overline{D}(\sum\limits_{\mathbb{N}}W^{m}_{n})$. Repeat the above proof method , we get that $\overline{D}(K_{m}(Z))=\overline{KF}(K_{m}(Z))$. Hence, $$\overline{D}(K_{m}(Z))=\overline{K}(K_{m}(Z))=\overline{KF}(K_{m}(Z))\mbox{ and }D_{m+1}(Z)=K_{m+1}(Z)=KF_{m+1}(Z).$$

For (3), let $m=\alpha$. By (2), we deduce $$D_{\alpha+1}(Z)=K_{\alpha+1}(Z)=KF_{\alpha+1}(Z).$$
Since $Z$ is $(\alpha+1)^{d}$-special and $(\alpha+1)^{wf}$-special, we have that for any ordinal $\delta<\alpha+1$, $K_{\delta}(Z)$ is not a $d$-space and $K_{\alpha+1}(Z)$ is well-filtered. Thus $K_{\delta}(Z)$ is not $k$-well-filtered and $K_{\alpha+1}(Z)$ is $k$-well-filtered. Then $Z$ is $(\alpha+1)^{k}$-special.

Next, we analyze the space $Z^{\top}$. For $Z^{\top}$, the statement (1) also follows from Lemma \ref{d rank} (2) and Lemma \ref{wf rank 2} (2).

For (2), The proof is by induction on $m$.

Basic steps. For $m=0$, clearly, $\overline{D}(Z^{\top})\subseteq \overline{KF}(Z^{\top})$. Conversely, let $A\in \overline{KF}(Z^{\top})$. From Lemma \ref{Liu bei}, we have that $A\in \overline{KF}(Z)$ or $A\sim \{\top\}$. It is straightforward to check that $A\in \overline{D}(Z^{\top})$.  Therefore, $$\overline{D}(Z^{\top})= \overline{KF}(Z^{\top})\mbox{ and }D(Z)=K(Z)=KF(Z).$$

Inductive steps. There are two cases to consider:

Case 1. Let $m$ be an ordinal such that $m+1<\alpha$. Assume that $$\overline{D}(K_{m}(Z^{\top}))=\overline{KF}(K_{m}(Z^{\top}))\mbox{ and }D_{m+1}(Z^{\top})=K_{m+1}(Z^{\top})=KF_{m+1}(Z^{\top}).$$

By the proof of Lemma \ref{d rank} (2), we have that $D_{m+1}(Z^{\top})=K_{m+1}(Z^{\top})\cong
 \sum\limits_{\mathbb{N}^{\top}}W^{m+1}_{n^{\prime}}$, where \[W^{m+1}_{n^{\prime}}=
 \begin{cases}
 D_{m+1}(X_{n^{\prime}}),& \mbox{if } m+1<\overline{\alpha}+n^{\prime}<\alpha,\\

 \mathrm{H}_{d}(X_{n^{\prime}}),& \mbox{if } \overline{\alpha}+n^{\prime}\leq m+1<\alpha, \\

 \mathbb{\top},& \mbox{if } n^{\prime}=\top. \\

\end{cases}\]
To prove $\overline{D}(K_{m+1}(Z^{\top}))=\overline{KF}(K_{m+1}(Z^{\top}))$, it is sufficient to show  $\overline{D}(\sum\limits_{\mathbb{N}^{\top}}W^{m+1}_{n^{\prime}})=
\overline{KF}(\sum\limits_{\mathbb{N}^{\top}}W^{m+1}_{n^{\prime}})$. Clearly, $\overline{D}(\sum\limits_{\mathbb{N}^{\top}}W^{m+1}_{n^{\prime}})\subseteq
\overline{KF}(\sum\limits_{\mathbb{N}^{\top}}W^{m+1}_{n^{\prime}})$. For any $A\in \overline{KF}(\sum\limits_{\mathbb{N}^{\top}}W^{m+1}_{n^{\prime}})$, again by Lemma \ref{Liu bei}, we have that$A\in\overline{D}(\sum\limits_{\mathbb{N}^{\top}}W^{m+1}_{n^{\prime}})$. Therefore, $\overline{KF}(\sum\limits_{\mathbb{N}^{\top}}W^{m+1}_{n^{\prime}})\subseteq
\overline{D}(\sum\limits_{\mathbb{N}^{\top}}W^{m+1}_{n^{\prime}})$. So $$\overline{D}(K_{m+1}(Z^{\top}))=\overline{KF}(K_{m+1}(Z^{\top})).$$
This implies
$$D_{m+2}(Z^{\top})=K_{m+2}(Z^{\top})=KF_{m+2}(Z^{\top}).$$

Case 2. Suppose that $m<\alpha$ is a limit ordinal and that the required statement holds for any $\delta<m$. Then
$$K_{m}(Z^{\top})=\bigcup\limits_{\delta<m}K_{\delta}(Z^{\top})\cong
\bigcup\limits_{\delta<m}D_{\delta}(Z^{\top})\cong
\bigcup\limits_{\delta<m}\sum\limits_{\mathbb{N}^{\top}}W^{\delta}_{n^{\prime}}\cong
\sum\limits_{\mathbb{N}^{\top}}\bigcup\limits_{\delta<m}W^{\delta}_{n^{\prime}}\cong
\sum\limits_{\mathbb{N}^{\top}}W^{m}_{n^{\prime}}.$$

Now it is enough to show $\overline{KF}(\sum\limits_{\mathbb{N}^{\top}}W^{m}_{n^{\prime}})=
\overline{D}(\sum\limits_{\mathbb{N}^{\top}}W^{m}_{n^{\prime}})$. Repeat the proof method of Case 1, we get that $\overline{D}(K_{m}(Z^{\top}))=\overline{KF}(K_{m}(Z^{\top}))$. So $$\overline{D}(K_{m}(Z^{\top}))=\overline{K}(K_{m}(Z^{\top}))=
\overline{KF}(K_{m}(Z^{\top})).$$
Hence
$$D_{m+1}(Z^{\top})=K_{m+1}(Z^{\top})=KF_{m+1}(Z^{\top}).$$

For (3), by (2), we get that
$$K_{\alpha}(Z^{\top})=\bigcup\limits_{\delta<\alpha}K_{\delta}(Z^{\top})\cong
\bigcup\limits_{\delta<\alpha}KF_{\delta}(Z^{\top})\cong KF_{\alpha}(Z^{\top}).$$
So $K_{\alpha}(Z^{\top})$ is well-filtered by (1) and hence $K_{\alpha}(Z^{\top})$ is $k$-well-filtered. Again by (1), we have that $K_{\delta}(Z^{\top})$ is not $k$-well-filtered for any ordinal $\delta<\alpha$. Therefore, $rank_{k}(Z^{\top})=\alpha$.
\end{proof}

\begin{theorem}\label{special exist}
For any ordinal $\alpha$, there exists a $T_{0}$ space $X$ and the following statements hold.
\begin{enumerate}
\item[(\texttt{1})] $X$ is $(\alpha+1)^{d}$-special and $(\alpha+1)^{wf}$-special.
\item[(\texttt{2})] $\overline{D}(K_{m}(X))=\overline{KF}(K_{m}(X))$ and $D_{m+1}(X)=K_{m+1}(X)=KF_{m+1}(X)$ for $0\leq m<\alpha+1$.
\item[(\texttt{3})] $X$ is $(\alpha+1)^{k}$-special.
\end{enumerate}
\end{theorem}

\begin{proof}
As usual, we use induction on $\alpha$. Base cases. If $\alpha=0$, by Lemma \ref{N is 1}, $(N, \tau_{\sigma})$ satisfies the required statements.

Inductive steps. There are two cases to consider:

Case (1). Suppose that the statements of the theorem hold for any ordinal $\alpha$. From Lemma~\ref{not a limit ordinal}, we have that the statements of the theorem hold for the ordinal $\alpha+1$.

Case (2). Let $\alpha$ be a limit ordinal. Assume for any ordinal $\beta<\alpha$, the statements of the theorem hold. Note that for every $n\in N$, $\overline{\alpha}+n<\alpha$. So the statements of the theorem hold for $\overline{\alpha}+n$. By Lemma \ref{limit rank}, it is straightforward to check that the statements of the theorem hold for $\alpha$.
\end{proof}

\begin{corollary}\label{non-limit k}
For any non-limit ordinal $\alpha$, there exists a $T_{0}$ space $X$ satisfying the following conditions.
\begin{enumerate}
\item[(\texttt{1})] $X$ is $\alpha^{d}$-special and $\alpha^{wf}$-special.
\item[(\texttt{2})] $\overline{D}(K_{m}(X))=\overline{KF}(K_{m}(X))$ and $D_{m+1}(X)=K_{m+1}(X)=KF_{m+1}(X)$ for $0\leq m<\alpha$.
\item[(\texttt{3})] $X$ is $\alpha^{k}$-special.
\end{enumerate}
\end{corollary}
\begin{proof}
For any non-limit ordinal $\alpha$, there exists an ordinal $\beta$ such that $\alpha=\beta+1$. By Theorem~\ref{special exist}, for ordinal $\beta$, there exists a $T_{0}$ space $X$ which satisfies the the required statements.
\end{proof}

\begin{theorem}\label{irreducible X}
For any ordinal $\alpha$, there exists an irreducible $T_{0}$ space $X$ such that  the following statements hold.
\begin{enumerate}
\item[(\texttt{1})] $rank_{d}(X)=rank_{wf}(X)=\alpha$.
\item[(\texttt{2})] $\overline{D}(K_{m}(X))=\overline{KF}(K_{m}(X))$ and $D_{m+1}(X)=K_{m+1}(X)=KF_{m+1}(X)$ for $0\leq m<\alpha$.
\item[(\texttt{3})] $rank_{k}(X)=\alpha$.
\end{enumerate}
\end{theorem}

\begin{proof}
Let $\alpha$ be an ordinal. Now we consider two cases:

Case 1. If $\alpha$ is not a limit ordinal, the statements of the theorem follow directly from Corollary~\ref{non-limit k}.

 Case 2. If $\alpha$ is a limit ordinal, then $\overline{\alpha}+n$ is a non-limit ordinal for any $n\in N$. Hence, the statements of the theorem follow from Lemma \ref{limit rank} and Corollary \ref{non-limit k}.
\end{proof}

\begin{corollary}\label{5}
For any ordinal $\alpha$, there exists an irreducible $T_{0}$ space $X$ whose $k$-rank is equal to $\alpha$.
\end{corollary}
\begin{proof}
For any ordinal $\alpha$, it is enough to show that the $k$-rank of the space $X$ required in Theorem~\ref{irreducible X} is $\alpha$.
\end{proof}

We can directly obtain the following corollaries.
\begin{corollary}(see also in \cite{Ershov 17})
For any ordinal $\alpha$, there exists an irreducible $T_{0}$ space $X$ whose $d$-rank is equal to $\alpha$.
\end{corollary}

\begin{corollary}(see also in \cite{Liu21})
For any ordinal $\alpha$, there exists an irreducible $T_{0}$ space $X$ whose $wf$-rank is equal to $\alpha$.
\end{corollary}

\bibliographystyle{plain}

\end{document}